\documentclass[12pt, a4paper]{article}
\usepackage{amsmath, amsthm,  amsfonts, amscd, epsfig, url}

\renewcommand{\a}{\alpha}
\renewcommand{\b}{\beta}
\renewcommand{\c}{\gamma}

\renewcommand{\P}{{\cal P}}

\newcommand{\Z}{\mathbb Z}

\newcommand{\C}{\mathbb C}

\newcommand{\G}{\mathcal G}

\newcommand{\cstar}{\C^\times}
\newcommand{\Cstar}{\C^\times}

\newcommand{\deligneone}{\underline{\cstar} \overset{d\log}{\to} \Omega^1 }
\newcommand{\delignetwo}{\underline{\cstar} \overset{d\log}{\to} \Omega^{1} \overset{d}{\to}\Omega^{2}}

\theoremstyle{plain}
\newtheorem{theorem}{Theorem}[section]

\newtheorem{proposition}[theorem]{Proposition}

\theoremstyle{definition}
\newtheorem{definition}[theorem]{Definition}

\theoremstyle{remark}

\newtheorem{note}{Note}[section]

\title{Bundle gerbes: stable isomorphism and local theory.}
\author{Michael K. Murray\thanks{The support of the Australian 
Research Council is acknowledged}\\
Department of Pure Mathematics\\
University of Adelaide\\
Adelaide, SA 5005 \\
Australia\\
\url{mmurray@maths.adelaide.edu.au}
\and
Daniel Stevenson\thanks{The support of an Australian Postrgraduate 
Award is  acknowledged}\\
Department of Pure Mathematics\\
University of Adelaide\\
Adelaide, SA 5005 \\
Australia\\
\url{dstevenso@maths.adelaide.edu.au}
}

\begin{document}
\maketitle

\begin{abstract} We consider the notion of stable isomorphism of 
	bundle gerbes. It has the consequence that the stable isomorphism
	classes of bundle gerbes over a manifold $M$ are in bijective 
	correspondence with $H^3(M, \Z)$.  Stable isomorphism sheds light
	on the local theory of bundle gerbes and  enables us to develop 
	a classifying theory for bundle gerbes using results of Gajer on 
	$B\Cstar$ bundles. 

\noindent AMS Classification scheme numbers: 55R65, 55R35

\end{abstract}

\section{Introduction}
Recently there has been interest in applying the theory of 
gerbes to differential geometry \cite{Hit} and physics 
\cite{Kal}.   In \cite{Mur} the first named author described 
the basic structure of bundle gerbes.  In as much as a gerbe is a 
sheaf of groupoids and a bundle gerbe is a bundle of groupoids, the
relationship gerbes and bundle gerbes  is analogous to that between 
locally invertible sheaves and line bundles. 

Just as  for  gerbes, a bundle gerbe over a manifold $M$ has a 
characteristic class, called its Dixmier-Douady class, which 
lives in $H^{3}(M, \Z)$ and vanishes precisely when the 
bundle gerbe is trivial.  In \cite{Mur} the theory of 
bundle gerbes was shown to be in many ways analogous to that
of complex line bundles with the Dixmier-Douady class playing
the role of the Chern class.  One missing link in 
 \cite{Mur} was the fact that  
the natural notion of isomorphism of bundle gerbes had the
unpleasant property that two bundle gerbes could be 
non-isomorphic but have the same Dixmier-Douady class. This, of 
course, is
not true for line bundles and their Chern class. In \cite{CarMicMur}
the notion of {\em stable isomorphism} of bundle gerbes was introduced.
This does have the property that two bundle gerbes are stably
isomorphic precisely when they have the same Dixmier-Douady
class.  In the present work we  give the precise definition 
of stable isomorphism and elaborate its properties.

A consequence of the definition of stable isomorphism is a local 
theory of bundle gerbes analogous to describing a line bundle
by transition functions. This description of bundle gerbes is 
the same as the local theory of gerbes used in \cite{Hit}.

\section{Bundle gerbes}
\label{bundle_gerbes}
Before recalling the definition of bundle gerbe from \cite{Mur} we
need some notation for fibre products and some constructions for
$\cstar$ bundles.  Mostly we will be working in the category of smooth 
manifolds and maps but often these will need to be 
infinite-dimensional manifolds and in the classifying theory they 
will be  differentiable spaces not manifolds \cite{Gaj}. 
In the interests of brevity we will just say map and assume the 
category can be deduced from context. 

We will be interested in maps  $\pi \colon Y \to M$ which admit local
sections.  That is, for every $x \in M$ there is an open set $U  $
containing $x$ and  a local section $s \colon U \to Y$. 
 For want of a better term we will call maps like this 
locally split.   Note that a locally split map is necessarily 
surjective and that if we are dealing with the smooth 
category a locally split map is just a submersion.
Locally trivial fibrations are, of course, locally split
but we will  not require $Y \to M$ to be a fibration as we did  
in \cite{Mur} because
there are a number of important examples in which $Y \to M$
only admits local sections. In particular if $U 
 = \{U_\a\}_{\a \in I}$ is an open cover over $M$
admitting local sections $s_\a \colon U_\a \to P$ let $Y_U$ be the 
disjoint union of all the elements in the open cover, that is the 
{\em nerve} of the open cover.  Then, of course, $Y_U \to M$ is 
rarely a fibration.

Let $Y \to M$ be locally split. Then we denote by
 $Y^{[2]} = Y\times_\pi Y$ the fibre product 
of $Y$ with itself over $\pi$, 
 that is the subset of pairs $(y, y')$ in
$Y \times Y$ such that $\pi(y) = \pi(y')$.  More generally we denote 
the $p$th fold fibre product by $Y^{[p]}$.

Let  $R$ and $S$ be principal $\Cstar$ spaces, that
is spaces on which $\cstar$ acts smoothly, freely and transitively.
Then there is a notion of a dual space $R^*$; this is the same space as 
$R$ but
with the inverse  $\Cstar$ action. The product  
 $R \times S$ is a principal $\Cstar \times \Cstar$ space.
If we quotient $R \times S$ by the subgroup $\{(z, z^{-1}) \mid z \in \Cstar \}$ 
we obtain the  contracted product $R \otimes S$ a principal $\Cstar$ 
space  \cite{Bry}. These  constructions extend in a fibrewise manner to principal
$\cstar$ bundles so that if $P$ and $Q$ are principal $\cstar$ bundles
we can form $P^*$ and $P \otimes Q$. 
If we replace $P$ and $Q$  by associated complex line bundles 
in the standard fashion then these two operations are just
linear dual and tensor product.

A bundle gerbe over $M$ is a pair $(P, Y)$ where
$\pi \colon Y \to M$ is a locally split map and 
$P$ is a $\cstar$ bundle $P \to Y^{[2]}$ with a product, that is,
a $\Cstar$ equivariant map
$$
P_{(y_1, y_2)} \otimes P_{(y_2, y_3)} \to P_{(y_1, y_3)}
$$
for every $(y_1, y_2)$ and $(y_2, y_3)$ in $Y^{[2]}$. 
We require the product to be smooth in $y_1$, $y_2$ and 
$y_3$ but in the interests of brevity we will not state the various
definitions needed to make this requirement precise, they  can be found in 
\cite{Mur}.
The product is required to be
associative whenever triple products are defined. We shall
often refer to a bundle gerbe $(P, Y)$ as just $P$. 
It is shown in \cite{Mur} that the product defines 
uniquely an inverse
$$
P_{(y_1, y_2)} \to P_{(y_2, y_1)}
$$
denoted by $p \mapsto p^{-1}$ and an identity element  $1 \in P_{(y, y)}$ which 
behave as one would expect. For example $(pq)^{-1}1 = q^{-1}p^{-1}$.

Various operations are possible on bundle gerbes.  Let $(P, Y)$ 
be a bundle gerbe over $M$. 
Let $\pi \colon Z \to N$ be another locally split map
and let  $\hat\phi \colon Z \to Y$ be a fibre map covering 
a map $\phi \colon N \to M$.  Then there is an induced map 
${\hat \phi}^{[2]} \colon Z^{[2]} \to Y^{[2]}$ which can 
 be used to pull-back the bundle $P \to Y^{[2]}$ to a bundle 
$({\hat\phi}^{[2]})^{-1}(P) \to Z^{[2]}$.  This has an induced product on it
and defines a bundle gerbe which we denote, for simplicity,
by $(\phi^{-1}(P), Z)$ or $\phi^{-1}(P)$.  Two special cases
of this are important. The first is when we just have just a map 
 $f\colon N \to M$ and use this to pull-back $Y \to M$ to 
$f^{-1}(Y) \to N$. The second is when we have $M = N$ and 
$\phi$ the identity. 

If $(P, Y)$ is a bundle gerbe  we can define a 
new bundle gerbe, $(P^*, Y)$, the dual of $(P, Y)$, by taking
the dual of $P$. 
Also if  $(P, Y)$ and $(Q, Z)$ are two bundle gerbes we can define their
product $(P\otimes Q,  Y\times_\pi Z)$ where $Y\times_\pi Z$ is 
the fibre product.

Two bundle gerbes $(P, Y)$ and $(Q, Z)$ are called {\em isomorphic} 
if there is an isomorphism  $f \colon Y \to Z$, commuting with the 
projection to $M$ and a bundle isomorphism $g \colon P \to Q$ covering the 
induced map $f^{[2]} \colon Y^{[2]} \to Z^{[2]}$ and commuting with the
bundle gerbe products on $Q$ and $P$ respectively.

If $Q$ is a $\cstar$ bundle over $Y$ then we can define a 
bundle gerbe $\delta(Q)$ by $\delta(Q) = {\pi_1^{-1}(Q)}^* \otimes 
\pi_2^{-1}(Q)$, that is $\delta(Q)_{(y_1, y_2)} = Q_{y_1}^* \otimes Q_{y_2}$. 
The bundle gerbe product 
is induced by the natural pairing
$$
Q_{y_1}^*\otimes Q_{y_2}\otimes Q_{y_2}^*\otimes Q_{y_3} \to 
Q_{y_1^*}\otimes Q_{y_3}.$$

A bundle gerbe which is isomorphic
to a  bundle gerbe of the form $\delta(Q)$ is  called {\em trivial}. 
A choice of $Q$ and  a bundle gerbe  isomorphism $\delta(Q) \simeq P$ is called
a {\em trivialisation}.  If $Q$ and $R$ are trivialisations
of $P$ then we have natural isomorphisms 
$$
Q_{y_1}^*\otimes Q_{y_2} \simeq R_{y_1}^*\otimes R_{y_2}
$$
and hence 
$$
Q_{y_1}^*\otimes R_{y_1} \simeq Q_{y_2}^*\otimes R_{y_2}
$$
so that the bundle $Q \otimes R$ is the pull-back of a
bundle on $M$. Moreover if $Q$ is a trivialisation and 
$S$ is a bundle on $M$ then $Q \otimes \pi^{-1}(S)$ is 
also a trivialisation.

There is a characteristic class $d(P) = d(P, Y) \in H^3(M, \Z)$, 
the Dixmier-Douady class of $(P, Y)$.  We recall the definition 
from \cite{Mur} as we will need it later.
Choose an open cover $\{ U_\alpha \}_{\a \in I}$ of 
$M$ for which all intersections are contractible.  Let $s_\a \colon 
U_\a \to Y$ be sections. These induce sections $(s_\a, s_\b) \colon 
U_\a \cap U_\b \to Y^{[2]}$. Choose sections $\sigma_{\a\b}$ of
$(s_\a, s_\b)^{-1}(P)$ over $U_\a \cap U_\a$.  The bundle gerbe 
multiplication allows us to compare $\sigma_{\a\b}\sigma_{\b\c}$ and
$\sigma_{\a\c}$ and define a map $g_{\a\b\c} \colon U_\a \cap U_\b 
\cap U_\c \to \C^\times$ by 
$$
\sigma_{\a\b}\sigma_{\b\c} =  \sigma_{\a\c}  g_{\a\b\c}.
$$
This map $g_{\a\b\c}$ is a 
$\C^\times$ valued cocycle and defines the Dixmier-Douady class in 
$H^2(M, \C^\times) \simeq H^3(M, \Z)$.

It is shown in \cite{Mur} that
\begin{theorem}[\cite{Mur}]
\label{th:trivial}
A bundle gerbe $(P, Y)$ has zero Dixmier-Douady class
precisely when it is trivial.
\end{theorem}

The construction 
of the Dixmier-Douady class is natural in the sense that 
if $Z \to N$ is another locally split map
and $\hat\phi \colon Z \to Y$ is a  fibre map covering 
$\phi \colon N \to M$ then it is straightforward to 
check from the definition that 
\begin{equation}
\label{eq:natural0}
d(\phi^{-1}(P), Z) = \phi^*(d(P, Y)).
\end{equation}
In particular if $M = N$ and $\phi $ is the identity 
then 
\begin{equation}
\label{eq:natural1}
d(\phi^{-1}(P)) = d(P).
\end{equation}

 From \cite{Mur} we also have
\begin{theorem}[\cite{Mur}]
\label{th:dd}
If $P$ and $Q$ are bundle gerbes over $M$ then
\begin{enumerate}
\item $d(P^*) = -d(P)$ and
\item$d(P\otimes Q) = d(P)
+ d(  Q)$. 
\end{enumerate}
\end{theorem}

\section{Stable isomorphism of bundle gerbes}
Equation \eqref{eq:natural1} shows that there are many 
bundle gerbes which have the same Dixmier-Douady class but 
which  are not isomorphic. 
We consider in this 
section the notion of {\em stable isomorphism}  of bundle gerbes ---
introduced  in  \cite{CarMicMur} --- that 
will have the consequence that two bundle gerbes have the same
Dixmier-Douady class if and only if they are stably isomorphic.
We will 
show in a later section that there is a {\em universal bundle gerbe}
with the property that every bundle gerbe over $M$
is stably isomorphic to the pull-back of the 
universal bundle gerbe under a classifying map.

We begin with 
\begin{definition}
\label{def:stableiso}
Two bundle gerbes $(P, Y)$ and $(Q, Z)$ are
called {\em stably isomorphic} if there are trivial bundle 
gerbes $T_1$ and $T_2$ such that 
$$
P \otimes T_1 = Q \otimes T_2.
$$
\end{definition}
There are equivalent definitions of stable isomorphism
provided by the following proposition.

\begin{proposition}
\label{prop:stable}
For bundle gerbes $(P, Y)$ and $(Q, Z)$
the following are equivalent.
\begin{enumerate}
\item $P$ and $Q$ are stably isomorphic
\item $P \otimes Q^*$ is trivial
\item $d(P) = d(Q)$. 
\end{enumerate}
\end{proposition}
\begin{proof}
Clearly stably isomorphic bundle gerbes have the same
Dixmier-Douady class because trivial bundles have the 
zero Dixmier-Douady class and the Dixmier-Douady class is 
additive over tensor products. So (1) implies (3).   If $d(P) = d(Q)$ then
$d(P\otimes Q^*) = d(P) - d(Q) = 0$. 
Hence $P\otimes Q^*$ is trivial (\cite{Mur}). So (3) implies (2). 
Finally if  $P \otimes Q^*$ is trivial then 
 $Q \otimes Q^*$ is also trivial as it has zero Dixmier-Douady class
and  then 
$P \otimes (Q^* \otimes Q) = Q \otimes (P \otimes Q^*)$
so $P$ and $Q$ are stably isomorphic. So (2) implies (1).
\end{proof}

Note that stable isomorphism for bundle 
gerbes is analogous to stable isomorphism
for vector bundles \cite{Ati}.

From part (3) of Proposition \ref{prop:stable} we see that stable 
isomorphism is an equivalence relation. 
It was shown in \cite{Mur} that every class in $H^3(M, \Z)$ is a
the Dixmier-Douady class of some bundle gerbe. Hence we can deduce
from Proposition \ref{prop:stable} that
\begin{theorem}
\label{th:stableiso}
The Dixmier-Douady class defines a bijection between
stable isomorphism classes of  bundle gerbes and $H^3(M, \Z)$. 
\end{theorem}

The following is an important example of stably
isomorphic bundle gerbes. 

\begin{proposition}
Let $\pi \colon Z \to M$ and $\pi \colon Y \to M$ be locally split 
maps and
 $\phi \colon Z \to Y$  a fibre map covering the identity on $M$. Let $(P, Y)$ be a bundle
gerbe then $(\phi^{-1}(P), Z)$ and
$(P, Y)$ are stably isomorphic. 
\label{prop:pullback}
\end{proposition}
\begin{proof} This is clear because they have the same
Dixmier-Douady class but it is instructive to also show that $\phi^{-1}(P)^* \otimes P$
is trivial. 
To see this note that
the fibre of $(\phi^{-1}(P))^* \otimes P$ at a point
$((z_1, y_1), (z_2, y_2)) $ in $(Z \times_\pi Y)^{[2]}$
is 
\begin{equation*}
\label{eq:pullback1}
\phi^{-1}(P)_{(z_1, z_2)}^*  \otimes P_{( y_1, y_2)} = 
P_{(\phi(z_1), \phi(z_2))}^*  \otimes P_{( y_1, y_2)}.
\end{equation*}
Whereas if we define $Q$ over $Z \times_\pi Y$ by 
$Q_{(z, y)} = P_{(\phi(z), y)}$ then 
$\delta(Q)$ at $((z_1, y_1), (z_2, y_2)) $ is
\begin{equation*}
\label{eq:pullback2}
P_{(\phi(z_1), y_1)}^*  \otimes P_{( \phi(z_2), y_2)}.
\end{equation*}
The gerbe multiplication can be used to define an isomorphism
\begin{align*}
P_{(\phi(z_1), y_1)}^*  \otimes P_{( \phi(z_2), y_2)} &=
P_{(\phi(z_1), \phi(z_2))}^* \otimes P_{(\phi(z_2), y_1)}^* \otimes P_{( \phi(z_2), 
y_2)}\\
&= P_{(\phi(z_1), \phi(z_2))}^* \otimes P_{(y_1, \phi(z_2))} \otimes P_{( \phi(z_2), 
y_2)}\\
&= P_{(\phi(z_1), \phi(z_2))}^* \otimes P_{(y_1, y_2)} 
\end{align*}
\end{proof}

We remark that an alternative approach to stable isomorphism is to 
define a stable isomorphism between bundle gerbes $P$ and $Q$
to be a trivialisation of $P \otimes Q^*$. Two bundle gerbes are
then stably isomorphic if a stable isomorphism exists. With a little
work one can show that stable isomorphisms are composable and that the set of 
all bundle gerbes over $M$ with stable isomorphisms as morphisms is a category,
in fact a two-category \cite{Ste}. We shall  not pursue this approach
here.  

\subsection{Local bundle gerbes}
Proposition \ref{prop:pullback} gives us a tool
for understanding what a bundle gerbe looks like locally. 
Consider a bundle gerbe $(P, Y)$ over $M$ and
assume that  $U = \{U_\a\}_{\a \in I}$ is an open cover over $M$
with local sections $s_\a \colon U_\a \to Y$.  Let $Y_U$ be the 
disjoint union of all the elements in the open cover, that is the 
{\em nerve} of the open cover.  Then the local sections $s_\a$ define 
a map $s \colon Y_U \to Y$ by $s(\alpha, x) = s_\alpha(x)$.  By
 \ref{prop:pullback} the 
pullback of the bundle gerbe over $Y$ defines a bundle gerbe over 
$Y_U$ which is stably isomorphic to the original $P$. 
If we specialise the definition of bundle gerbes to 
nerves of open covers we see that it is equivalent to the
following:
\begin{itemize}
\item   a $\cstar$ bundle $P_{\a\b}$ over each intersection $U_\a  \cap U_\b$
\item a trivialisation $\theta_{\a\b\c}$ of the contracted product 
     $P_{\b\c}\otimes P_{\a\c}^* \otimes P_{\a\b}$ over $U_\a \cap U_\b \cap U_\c$
\item the trivialisation satisfies $\delta(\theta) = 1$ over
$U_\a \cap U_\b \cap U_\c \cap U_\delta$.
\end{itemize}
Note in the last of these that  $\delta(\theta)$ is a section 
of a $\cstar$ bundle which is canonically trivial so the statement
makes sense.   This data can be 
easily  shown to be equivalent to that in \cite{Hit}.   The discussion we have
just given shows that every bundle gerbe is stably isomorphic
to a `local' bundle gerbe of this form.

\subsection{Other related objects}
The geometric structure underlying a bundle gerbe arises
in a number of other places which we briefly mention here.

Firstly we remark that it is well-known that 
the collection of spaces $Y^{[p]}$ form a {\em simplicial space} 
\cite{Dup}.  In \cite {BryMac}  Brylinski and Maclaughlin  define the notion of 
a {\em simplicial line bundle}. For the simplicial space $Y^{[p]}$
this is precisely the notion of a bundle gerbe. 

Secondly  we can regard  $Y$ as the set of objects for 
a non-transitive $\Cstar$ groupoid in which two objects
$y_1$ and $y_2$ either have no morphism between them
if $y_1$ and $y_2$ are not in the same fibre or
the single morphism $(y_1, y_2)$ if $y_1$ and $y_2$
are in the same fibre.   Composition of morphisms 
is $(y_1, y_2) \circ (y_2, y_3) = (y_1, y_3)$.  The
bundle gerbe $P$ is then a  central extension 
of the groupoid $Y$ called a {\em twist} \cite{Kum}.

Thirdly in a similar manner to the above any fibre $Y_m = \pi^{-1}(m)$ can
be thought of as a groupoid with the morphisms between $(y_1, y_2)$
being just $P_{(y_1, y_2)}$.  The bundle gerbe can then be thought of as a 
bundle of groupoids. 

\section{Bundle gerbes with connection}
In \cite{Mur} the first author defined  the notion 
of a bundle gerbe  connection, curving and curvature which 
we breifly recall. 

Recall from \cite{Mur} that if $\Omega^p(X)$ denotes the space
of $p$ forms on the manifold $X$, then  we have an 
{\em exact } complex
\begin{equation}
\label{eq:exact}
\Omega^p(M) \overset{\pi^*} {\to} \Omega^p(Y) 
\overset{\delta}{\to} \Omega^p(Y^{[2]}) \overset{\delta}{\to} \dots.
\end{equation}
Here $\delta \colon \Omega^p(Y^{[q]}) \to \Omega^p({Y^{[q+1]}})$
is the alternating sum of pull-backs $\sum_{j=1}^{q+1} (-1)^j 
\pi_i^*$ of projections where $\pi_i$ is the projection map which
omits the $i$th point in the fibre product. 
In \cite{Mur} exactness of the complex  \eqref{eq:exact} was proved for the case of $Y \to M$ being
a fibration but that proof can be extended to the case that $Y \to M$
is locally split.

Because $P \to Y^{[2]}$ is a $\Cstar$ bundle
it has connections. It is shown in \cite{Mur} that it
admits {\em bundle gerbe connections}
that is connections commuting  with the bundle gerbe product.
A bundle gerbe connection $\nabla$ has curvature $F_\nabla$
satisfying $\delta(F_\nabla) = 0 $ and hence
from  the exactness of equation \eqref{eq:exact} we have  that 
there exists a two-form $f$ on $Y$, satisfying 
the `descent equation'
$$
F_\nabla =  \pi_1^*(f) - \pi^*_2(f).
$$
Such an $f$ is called a {\em curving} for the connection $\nabla$.
The choice of a curving is not unique, from \eqref{eq:exact} we 
see that the ambiguity in the choice is precisely the addition of 
the pull-back of a two-form from $M$.
Given a choice of curving we then have that
 $\delta(df) = d \delta(f) = dF_\nabla = 0$
so that we can find some $\omega$, a three-form on $M$,
such that $df =  \pi^*(\omega)$. Moreover $\omega$ is closed as 
$\pi^*(d\omega) = ddf = 0$.  In \cite{Mur} it is shown that
$\omega/{2\pi i}$ is a de Rham representative for
the Dixmier-Douady class.  We call $\omega$ the 
three curvature of the connection and curving.

To define the notion of stable equivalence of
bundle gerbes with connection and curving we start
by letting   $P \to Y$ be a line
bundle with connection $\nabla$ and curvature $F$.
The trivial gerbe $\delta(P)  \to Y^{[2]}$
 has a natural  bundle gerbe connection $\delta(\nabla) = \pi_1^*(\nabla) -
 \pi_2^*(\nabla)$ and curving $F$.  If $T$ is a trivial bundle gerbe 
 we say it has trivial connection and curving if they arise in this 
 way.  To extend  the definition
of stable isomorphism (definition \ref{def:stableiso}) to cover the case of
bundle gerbes with connection and curving we  assume that the trivial bundle
gerbes $T_1$ and $T_2$, in the definition,  have  trivial connections 
and curvings and that the
isomorphism  in Definition \ref{def:stableiso} preserves connections 
and curving. 

\subsection{Bundle gerbes and Deligne cohomology}
In \cite{Mur} it was shown that a bundle gerbe with connection and curving defined
an element in the Deligne cohomology group $H^2(M, \delignetwo)$. In this
section we use the notion of stable isomorphism to show 
\begin{theorem}
\label{th:stabledeligne}
 The stable isomorphism classes of bundle gerbes with connection and
curving are equal to $H^2(M, \delignetwo)$.
\end{theorem}

Let us first define the Deligne cohomology class $D(P, \nabla, f)$ of a bundle gerbe 
$P$ with  connection $\nabla$ and curving $f$. 
If  ${\mathcal U} = \{ U_\a \}$ is a Leray cover of $M$ then the 
 Deligne cohomology $H^2(M, \delignetwo)$ can be calculated 
as the cohomology of the total complex of the double complex:
\begin{equation}
\label{double_complex}
\begin{array}{ ccccc}
\ \phantom{\delta} \vdots &  &  \ \phantom{\delta}\vdots && \ \phantom{\delta}\vdots
\\
 \delta\uparrow &  & \delta\uparrow &&\delta\uparrow 
\\
C^2({\mathcal U}, \underline{\cstar}) &  \stackrel{d\log}{\to} &
C^2({\mathcal U}, \Omega^1) & \stackrel{d }{\to} & C^2({\mathcal U}, \Omega^2) 
 \\
 \delta\uparrow &  & \delta\uparrow & &\delta\uparrow
\\
C^1({\mathcal U}, \underline{\cstar}) &  \stackrel{d\log}{\to} &
C^1({\mathcal U}, \Omega^1) & \stackrel{d}{\to} & C^1({\mathcal U}, \Omega^2) 
 \\
 \delta\uparrow &  & \delta\uparrow & &\delta\uparrow \\
 C^0({\mathcal U}, \underline{\cstar}) &  \stackrel{d\log}{\to} &
C^0({\mathcal U}, \Omega^1) & \stackrel{d }{\to} & C^0({\mathcal U}, \Omega^2) 
 \\
\end{array}
\end{equation}

Consider an open cover ${\mathcal U} = \{ U_\a \}$ of $M$ which has 
all intersections contractible and 
 such that $Y \to P$ has sections $s_\a$ over each $U_\a$. 
Then $(s_\a, s_\b )$ defines a section of $Y^{[2]}$ over $U_\a \cap  U_\b$ 
which can be used to pull back $P$ to a bundle $P_{\a\b}$. Let 
$\sigma_{\a\b}$ be a section of this bundle. We have seen that the Dixmier-Douady 
class is the class $g_{\a\b\c}$ 
defined by 
$\sigma_{\a\b}\sigma_{\b\c} = g_{\a\b\c} \sigma_{\a\c}$.
  The connection $\nabla$ can be pulled back  
by $\sigma_{\a\b}$ to 
define $A_{\a\b}$ a 1-form on $U_\a \cap U_\b$ and the curving can be
pulled-back by the $s_\a$ to define two-forms $f_\a $ on $U_\a$. The 
triple $(g_{\a\b\c}, A_{\a\b}, f_\a)$ is a Deligne cohomology class 
which  we denote by $D(P, \nabla, f)$.  It is straightforward to check
that it is independent of choices.  When it is clear what the 
connection and curving are we shall often denote the Deligne
class by just $D(P)$. 

To prove Theorem \ref{th:stabledeligne} we mimic the proof of Theorem 
\ref{prop:stable}. 
Let $D(P)$ denote the Deligne class of a bundle gerbe with connection 
and curving.  
Note first that it follows easily from the constructions in \cite{Mur}
that  if $P$ and $Q$ are bundle gerbes with connection and 
curving and we give $P \otimes Q$ and $P^*$ the induced connections
and curving then $D(P \otimes Q) = D(P) + D(Q)$ and $D(P^*) = -D(P)$. 
Assume that we have:

\begin{proposition}
\label{prop:trivial}
Let $P$ be a bundle gerbe with connection and curving and zero 
Deligne class. Then $P$ is a trivial bundle gerbe with trivial 
connection and curving. 
\end{proposition}

Consider now  $P$ and $Q$  two bundle gerbes with connection and 
curving and equal Deligne classes then we can form $P \otimes Q^*$. 
We have $D(P \otimes Q^*) = 0$ so by  Proposition \ref{prop:trivial} 
this is a trivial bundle gerbe with trivial connection and curving. 
Similarly for $Q \otimes Q^*$. But then we have 
$P \otimes (Q^* \otimes Q) =  Q \otimes (P \otimes Q^*) $ so that 
$P$ and $Q$ are stably isomorphic. 

It remains to prove Proposition \ref{prop:trivial}. 
Assume that  $(P, Y)$ is a bundle gerbe over $M$ with 
$D(P)= 0$.  This implies that the Dixmier-Douady class
of $P$ is also zero and we may assume that $P = \delta(Q)$ for
$Q \to Y$ a $\cstar$ bundle.   Let $\nabla_P$ be the 
connection on $P$ and $\nabla_Q$ a connection on $Q$. Then 
$\delta(\nabla_Q)$ is another connection on $P$ and hence
$\nabla_P - \delta(\nabla_Q)$ is a 1-form with 
$\delta(\nabla_P - \delta(\nabla_Q)) = 0$. By the exactness of 
\ref{eq:exact} there is 
a 1-form $a$ on $Y$ such that $\delta(a) = \nabla_P - \delta(\nabla_Q)$
or $\nabla_P = \delta(\nabla_Q - a)$.  Replacing $\nabla_Q$ by 
$\nabla_Q - a$ we may assume,  without 
loss of generality that $\nabla_P =\delta(\nabla_Q)$.  Consider now 
the  curving $f$ and the curvature $F_Q$ of $\nabla_Q$.  We have
$\delta(f - F_Q) = 0$ and hence $f - F_Q  = \pi^*\phi$ for $\phi $ a 
2-form  on $M$. 

Consider the construction of the Deligne class $D(P, \nabla_P, f)$. If we choose
sections $\sigma_\a $ of $s_\a^{-1}(Q)$ and let 
$\sigma_{\a\b} = \sigma_\a/ \sigma_\b$ we have $g_{\a\b\c} = 1$. 
If $a_\a$ is the connection 1-form for $s_\a^{-1}(\nabla_Q)$
then $A_{\a\b} = a_\a - a_\b$. The vanishing of the 
Deligne class implies the existence of 
 $h_{\a\b}:U_{\a\b} \to \cstar$ 
and $k_{\a} \in \Omega^{1}(U_{\a})$ such that 
\begin{eqnarray*}
g_{\a\b\c}&  = & h_{\b\c} h_{\a\c}^{-1} h_{\a\b}        \\
A_{\a\b}&  =&  -k_{\b}+k_{\a}+h_{\a\b}^{-1}dh_{\a\b}  \\
f_{\a} & = & dk_{\a}                                \\
\end{eqnarray*}
and hence in our case
\begin{eqnarray*}
1&  = & h_{\b\c} h_{\a\c}^{-1} h_{\a\b}       \\
a_\a - a_\b &  =&  -k_{\b}+k_{\a}+h_{\a\b}^{-1}dh_{\a\b}  \\
f_{\a} & = & dk_{\a}.                              \\
\end{eqnarray*}

We can use the $h_{\a\b}$ to construct a bundle $R \to M$ with 
a connection $\nabla_R$ defined locally by $ k_\a - a_\a $. 
Pull this back to $Y$ and tensor it by $Q$. Then 
 $\delta(Q \otimes \pi^{-1}(R)) = \delta(Q) = P$. 
Moreover the connection $k_\a - a_\a$ on $R$ has
curvature $s_\a^*(f_\a - F_Q) = s_\a^* \pi^*(\phi) = \phi$. 
Hence the curvature of $\pi^*(\nabla_R)$ is $\pi^*(\phi) = f - F_Q$ and
when we tensor $\nabla_R$ with $\nabla_Q$ the resulting tensor 
product  connection 
will have curvature $f$.  We have now shown that a bundle
gerbe with trivial Deligne cohomology class is a trivial bundle
gerbe with trivial connection and curving.  

\begin{note}
Note that in the proof of Proposition \ref{prop:trivial} we have shown that if $P = \delta(Q)$ is a 
trivial bundle
gerbe with connection $\nabla$ then 
there is a connection on $Q$ such that $\delta(Q) = \nabla$. This 
implies that stable isomorphism classes of bundle gerbes with 
connection (but not curving) are in bijective 
correspondence with  stable isomorphism 
classes of bundle gerbes. As the former are determined by their 
Deligne class in  $H^2(M, \deligneone)$ and the latter by their 
Dixmier-Douady class in  $ H^2(M, 
\cstar)$ this is equivalent to the easily proven result that $H^2(M, \deligneone) = 
H^2(M, \cstar)$. 
\end{note}

\section{Bundle gerbes and gerbes}
In \cite{Mur} the first author proposed a construction 
of a presheaf of categories from a bundle gerbe whose
sheafification, it was claimed, was a gerbe associated
to the bundle gerbe and with the same Dixmier-Douady class. 
We give here a simpler construction of a gerbe associated
to a bundle gerbe which includes the previous case as 
a subsheaf.  Let $(P, Y) $ be a 
gerbe over $M$.  For any $U$ open in $M$ let 
$\G(U) $ be a category defined as follows. 
The objects of $\G(U)$ are the set of all trivialisations of the 
restriction of $(P, Y)$ to $U$. That is all pairs $(Q, f)$
where $Q$ is a $\Cstar$ bundle over $Y_U = \pi^{-1}(U) \subset Y$
and $f \colon \delta(P) \to Q_{|Y_U^{[2]}}$ is an isomorphism 
of bundle gerbes.  The morphisms between two objects 
$(Q, f)$ and $(R, g)$ are all isomorphisms of bundle 
gerbes which commute with $f$ and $g$. That is isomorphisms 
$\phi\colon Q \to R$ of $\Cstar$ bundles, acting 
as the identity on $Y_U$ such that the 
induced map 
$$
\delta(\phi) \colon \delta(Q) \to \delta(R)
$$
commutes with the isomorphisms $f$ and $g$.

Notice that $\G(U)$ is either empty or 
a groupoid and from section \ref{bundle_gerbes} we see that 
 $\G(U)$ is non-empty if $Y $ admits a section over
 $U$. 
 
We want to show that $\G$ defines a gerbe. The 
definition is in \cite{Bry} except that we will 
not go for the full generality of assigning to 
any local homeomorphism a groupoid although that 
would not be difficult. 

 First we need to show
that it is a sheaf of groupoids. We have already
shown that for every open set $U$ we have a (possibly trivial)
groupoid. The restriction functor is just exactly the
restriction of a trivialisation over $Y_U$ to 
$Y_V$ if $V \subset U$. For a general sheaf of
categories it is possible that the composition 
of two restriction functors differs from the restriction
functor of the composition of two restrictions by a 
natural transformation. In the case of $\G$ they
are clearly equal. This makes $\G $ a presheaf of 
groupoids (\cite{Bry} Definition 5.2.1). 

To show that $\G$ is a sheaf  of groupoids we need to 
check two patching conditions on objects and morphisms 
as detailed on pages 191--2 of \cite{Bry}.
Assume we have an open cover $\{ U_\a\}$ of an 
open set $U$. First consider  two trivialisations $(Q_i, f_i)$, 
$i=1,2$ in $\G(U)$ with morphisms $\phi_\a \colon {Q_1}_{|U_\a} \to {Q_2}_{|U_\a} $
for each $\a$ agreeing on overlaps. Then these clearly 
patch together to yield a global morphism $\phi$ and as the
$\phi_\a$ commute with the $f_i$ so also does $\phi$.  Second 
assume we have trivialisations $(Q_\a, f_\a)$ in each $\G(U_\a)$
and morphisms
$$
\phi_{\a\b} \colon {Q_\a}_{|U_{\a\b}} \to {Q_\b}_{|U_{\a\b}}
$$
satisfying $\phi_{\a\b}\phi_{\b\c}\phi_{\c\a} = 1$. Then we need to 
find a global trivialisation $(Q, f) \in \G(U)$ whose
restriction to each $U_\a$ is $(Q_\a, f_\a)$. But this
is possible because the condition $\phi_{\a\b}\phi_{\b\c}\phi_{\c\a}= 1$
allows us to apply the `clutching construction' to form the 
bundle $Q$. It is straightforward to define the trivialisation.
Hence we have that $\G$ is a sheaf of categories. 

Finally we consider the conditions in \cite{Bry} 
Definition 5.2.4  that we need to satisfy to show that $\G$
is a gerbe. First condition (G1). Assume that $\G(U)$ is 
non-empty. Let $(Q, f)$ be an object in $\G(U)$ and 
consider the automorphisms of $(Q, f)$. If we think of 
$Q$ first as a $\Cstar$ bundle on $Y_U$ then 
the group of all automorphisms is the group of all
maps from $Y_U $ to $\Cstar$.  However if we require
that they also commute with $f$ it is easy to see that
they have to be maps that are constant on the fibres of
$Y \to M$. Hence they are the group of all maps from $U$
into $\Cstar$ as required to satisfy (G1). Second condition (G2).
Let $(Q, f)$ and $(R, g)$ be objects in $\G(U)$ and let $z \in U$.
We have that $Q \otimes R^* = \pi^{-1} (T)$ for some bundle $T$
over $U$. Choosing a contractible neighbourhood $V$ of 
$z$ we can trivialise $T$ and this induces an isomorphism 
from $Q_{|V}$ to $R_{|V}$ as required. 
Finally the third condition that we can cover $M$ by 
open sets $U$ such that $\G(U)$ is non-empty 
follows from the fact that we can cover by open sets 
over which $Y$ has sections and hence we can trivialise
the bundle gerbe locally.

To see that $\G$ has the same Dixmier-Douady class as $(P, Y)$
recall  that inside any $\G(U)$ is a subcategory 
whose objects are induced by choosing sections $s \colon U \to Y$.
  If we construct
the Dixmier-Douady class of the gerbe $\G$ by choosing objects
and morphisms of this kind we are constructing the 
Dixmier-Douady class of the bundle gerbe $(P, Y)$ \cite{Mur}.

Hence we have proven.

\begin{theorem}
The sheaf of groupoids $\G$ is a gerbe
with Dixmier-Douady class equal to that of $P$.
\end{theorem}

Because both stable isomorphism classes of bundle gerbes
and equivalence classes of gerbes are classified by $H^3(M, \Z)$
it follows that we have established a bijection between 
stable isomorphism classes of bundle gerbes and equivalence 
classes of gerbes.

\section{Bundle gerbes and $B\cstar$ bundles.}
In this section we recall Gajer's construction 
of $B\Cstar$ bundles \cite{Gaj} and show that there is a
bijection between $B\Cstar$ bundles over $M$ and 
stable isomorphism classes of bundle gerbes over $M$. 
Because both of these sets are in bijective correspondence
with $H^3(M, \Z)$ the only thing we really need to do 
is to associate to any $B\Cstar$ bundle a bundle gerbe 
whose Dixmier-Douady class is the characteristic
class  of the $B\Cstar$ bundle. This will also give
 us a classifying theory for bundle gerbes. To this end
we need to recall the notion of the {\em lifting bundle
gerbe}.

\subsection{Lifting bundle gerbes}
Let 
$$
0 \to \cstar \to \hat G \to G \to 0
$$
be a central extension of groups. If $P(M, G)$ 
is a principal bundle there is a map
$$
s \colon P^{[2]} \to G
$$
defined by $ps(p,q) = q$ for every $p$ and $q$ in the
same fibre of $P \to M$.  Denote by $\hat P$ the 
pullback of $\hat G$ by $s$. The group multiplication 
in $\hat G$ commutes with $\cstar \subset G$ as it is central
and makes $(\hat P, P)$ a bundle gerbe over $M$ called the 
{\em lifting bundle gerbe} of $P$. 
 It is shown 
in \cite{Mur} that the lifting bundle 
gerbe  is trivial if and
only if the bundle $P$ lifts to a $\hat G$ bundle. 
We shall see below that every gerbe is stably isomorphic
to a lifting bundle gerbe.

Note that the Chern class of $\hat G \to G$ 
defines a class in $H^2(G, \Z)$. If $G$ is connected
then this induces a unique class in $H^2(P_m, \Z)$
where $P_m = \pi^{-1}(m)$ for any $m$.  The transgression
map \cite{Spa} maps this to a class in $H^3(M, \Z)$.
We have from  \cite{CarCroMur}
\begin{proposition}
\label{prop:trans}
The transgression of the class in $H^2(P_m, \Z)$
defined by the Chern class of the bundle $\hat G \to G$
is the negative of the Dixmier-Douady class of the bundle gerbe 
$(\tilde P, P)$. 
\end{proposition}

\subsection{$B\Cstar$ bundles}
In \cite{Gaj} Gajer considers the groups $B^p\C^\times$
which are the $p$-fold classifying spaces of $\C^\times$.
To do this it is necessary to use Milgram's construction
of classifying spaces which has the property that 
if $G$ is an abelian group then $BG $ is an abelian
group. It follows that all the $B^p\C^\times$ are (abelian)
groups. It is also true that each $EB^{p-1}\C^\times$
is a group and there is a short exact sequence
$$
0 \to B^{p-1}\cstar \to EB^{p-1}\cstar \to B^p\cstar  \to 0
$$
of groups. 

Gajer shows that the isomorphism classes of $B^p\Cstar$ bundles
over a space $M$ are in bijective correspondence with the 
cohomology $H^{p+2}(M, \Z)$. In particular $H^3(M, \Z)$
is in bijective correspondence with the space of $B\Cstar$
bundles on $M$. We will show 
\begin{theorem}
The set of all stable isomorphism classes of bundle
gerbes is in bijective correspondence with the set of
all isomorphism classes of $B\Cstar$ bundles
on $M$
\end{theorem}
\begin{proof}
Because we have   the short exact sequence 
\begin{equation}
\label{eq:gajerseq}
0 \to \Cstar \to E\Cstar  \to B\Cstar \to 0
\end{equation}
we can associate to any $B\Cstar$ bundle $X$
over $M$ the equivalence class of its lifting bundle gerbe. We 
want to show that this 
map is a bijection. It is enough to show that the characteristic 
class of a $B\cstar$ bundle and the Dixmier-Douady class
of its associated lifting bundle gerbe are the same. 

The characteristic class of a $B\cstar$ bundle
is the pull-back of the generator of $BB\cstar$ 
by any classifying map.  Because $B^p\cstar$ is a $K(\Z, p+1)$
the transgression for the fibering $EB\cstar \to BB\cstar$
is an isomorphism $H^2(B\cstar, \Z) \to H^3(BB\cstar, \Z)$. 
Fix a  generator of $H^3(BB\cstar, \Z)$ by defining it to be 
the negative of the transgression of the Chern class of 
$E\cstar \to B\cstar$. 
A choice of classifying map defines a commutative
diagram
$$
\begin{array}{ccc}
X & \rightarrow  &EB\cstar \\
\downarrow & & \downarrow \\
M &\to& BB\cstar.
\end{array}
$$
It follows that transgression maps of these two fiberings
will commute with pull-backs and hence that the characteristic
class of the bundle $X \to M$ is the negative of the transgression of the 
generator of two-dimensional 
cohomology in any of the fibres. The required result follows
from \ref{prop:trans}.
\end{proof}

\section{Classifying theory for bundle gerbes}
\label{classifying}

The classifying theory for bundle gerbes now follows
from that for $B\Cstar$ bundles. We define the 
universal bundle gerbe to be the lifting bundle gerbe
for the fibration $EB\cstar \to BB\cstar$ and we have
\begin{theorem}
If $(P, Y)$ is a bundle gerbe over
$M$ there is a map $f \colon M \to BB\Cstar$, unique
up to homotopy, such that the pull-back of the 
universal bundle gerbe is stably isomorphic to $(P, Y)$.
\end{theorem}

Note that this theorem has as corollaries a number of results about 
the stable isomorphism classes of bundle gerbes.  For example it shows 
that  every  bundle gerbe is stably isomorphic to a lifting bundle gerbe 
and that  every 
bundle gerbe is stably isomorphic to a bundle gerbe $(P, Y)$ where the
map $Y \to M$ is a fibration with fibre $K(\Z, 2)$.

Given a class in $H^3(M, \Z)$ we can (following Gajer) use the exactness of 
\eqref{eq:gajerseq} to construct isomorphisms of sheaf cohomology 
groups:
$$
H^3(M, \Z) \simeq H^2(M, \underline{\cstar}) \simeq H^1(M, \underline{B\cstar})
\simeq H^0(M, \underline{BB\cstar})
$$
where the last of these spaces is the space of all maps from
$M$ to $BB\cstar$. Here we use the notation $\underline{G}$ for the 
sheaf of maps 
into the group $G$.  This gives a correspondence between
the Dixmier-Douady class of a bundle gerbe and the classifying 
map of the associated $B\cstar$ bundle. Explicit formulae for these maps, at the level
of cocycles,  are given in \cite{Gaj}
Remark 2.1. 

\subsection{Lifting bundle gerbes}
As an example consider the case of a lifting bundle gerbe. There we 
start with a principal $G$ bundle $P$ over $M$ where $G$ has a central
extension:
$$
0 \to \cstar \to \hat G \to G \to 0.
$$
To this we have seen there is associated a bundle gerbe and hence
a $B\Cstar$ bundle. On the other hand the classifying map 
for the $\cstar$ bundle $\hat G \to G$ is a map
$$
\chi  \colon G \to B\cstar.
$$
If $\chi$ was a group homomorphism then standard techniques 
would enable us to define an induced $B\cstar$ bundle
as $\P = P \times_G B\cstar$.  Such a $\chi$
exists if we work with geometric realisations 
of $BG$ instead of the fat realisations. 
This has the  advantage \cite{Seg}
that if $G$ is any group then $EG$ is a group
and $G \to EG$ is a homomorphism.  In particular
$\hat G \to E\hat G$ is a homomorphism and hence the 
induced map $G \to E\hat G / \cstar$ is a homomorphism. 
But because $E\hat G$ is contractible we must
have $E\hat G / \cstar$ homotopy equivalent to $B\cstar$.
For the purposes of the remaining discussion we take
$E\hat G / \cstar$ to be the realisation of $B\cstar$. 
We conclude that $B\cstar$ is a group and that
we have defined a homomorphism $\chi \colon G \to B\cstar$.
Because of the commuting diagram
$$
\begin{array}{ccc}
\hat G & \rightarrow & E{\hat G} \\
\downarrow &   &  \downarrow \\
G & \stackrel{\chi}{\rightarrow} & B\cstar 
\end{array}
$$
we see that the homomorphism $G \to B\cstar$ is a 
classifying map for the bundle $\hat G \to G$. 

The homomorphism $G \to B\cstar$ can be used to 
define an associated $\cstar$ bundle to $P$
by defining $\P = P \times_G B\cstar$. If we choose
a classifying map $\phi \colon M \to BG$ for
the bundle $P(M, G)$ then we have a commuting diagram
$$
\begin{array}{ccccc}
P        & \rightarrow &EG & \rightarrow & EB\cstar \\
\downarrow  &  &\downarrow   &  &\downarrow \\
M &\stackrel{\phi}{\rightarrow} & BG & \stackrel{B\chi}{\rightarrow} & 
BB\cstar
 \end{array}
$$

The composition $B\chi \circ \phi$ is a classifying map for the 
associated bundle $\P$ which is the $B\Cstar$ bundle associated
to the lifting bundle gerbe of $P$.

Note finally that if $G= LK$ is the group of smooth loops into a compact
Lie group then it is well-known that there is a commuting
diagram of $U(1)$ bundles
$$
\begin{array}{ccc}
\hat G & \rightarrow & U(\mathcal{F}) \\
\downarrow &   &  \downarrow \\
G & \rightarrow & PU\mathcal{F} 
\end{array}
$$
where $\hat G$ is the Kac-Moody central extension, $U(\mathcal{F})$
the unitary operators on the Fock space $\mathcal F$,
$PU(\mathcal{F})$ the projectivization of this group and the 
horizontal maps are all group homomorphisms \cite{PreSeg}.   Because 
$U(\mathcal{F})$ is contractible we see that $PU(\mathcal{F})$ is a 
realisation of $PU(1)$.  This gives a concrete realisation of the 
$BU(1)$ bundle associated to an $LK$  bundle.

\end{document}